\documentclass{amsart}
\usepackage{bm}
\usepackage{amsthm}

\newtheorem{theorem}{Theorem}[section]

\newtheorem{definition}[theorem]{Definition}

\newtheorem{lemma}[theorem]{Lemma}
\newtheorem{proposition}[theorem]{Proposition}

\newtheorem*{question}{Question}

\newcommand{\ilim}{\varprojlim}

\begin{document}

\title{Inverse Limits of Uniform Covering Maps}
\author{B.~LaBuz}
\address{Saint Francis University, Loretto, PA 15940}
\email{brendon.labuz@@gmail.com}
\date{December 22, 2009}

\begin{abstract}
In ``Rips complexes and covers in the uniform category'' the authors define, following James, covering maps of uniform spaces and introduce the concept of generalized uniform covering maps. Conditions for the existence of universal uniform covering maps and generalized uniform covering maps are given. This paper notes that the universal generalized uniform covering map is uniformly equivalent to the inverse limit of uniform covering maps and is therefore approximated by uniform covering maps. A characterization of generalized uniform covering maps that are approximated by uniform covering maps is provided as well as a characterization of generalized uniform covering maps that are uniformly equivalent to the inverse limit of uniform covering maps. Inverse limits of group actions that induce generalized uniform covering maps are also treated.
\end{abstract}

\maketitle
\tableofcontents

\medskip 
\medskip

\section{Introduction}

In \cite{Rips}, a theory of uniform covering maps and generalized uniform covering maps for uniform spaces is developed. In particular, it is shown that a locally uniform joinable chain connected space has a universal generalized uniform covering space and a path connected, uniformly path connected, and uniformly semilocally simply connected space has a universal uniform covering space. This paper points out that the universal generalized uniform covering map is uniformly equivalent to the inverse limit of uniform covering maps and is therefore approximated by uniform covering maps. Thus we consider when a generalized uniform covering map is approximated by uniform covering maps and when it is uniformly equivalent to the inverse limits of uniform covering maps. Inverse limits of  uniform covering maps in general are also investigated as are the inverse limit of group actions that induce uniform covering maps.

A good source for basic facts about uniform spaces is \cite{B}. Let us recall some definitions and results from \cite{Rips}. Given a function $f:X\to Y$ with $X$ a uniform space, the function \textbf{generates a uniform structure} on $Y$ if the family $\{f(E):E \mathrm{\ is\ an\ entourage\ of\ } X\}$ forms a basis for a uniform structure on $Y$. If $Y$ already has a uniform structure, the function generates that structure if and only if it is uniformly continuous and the image of every entourage of $X$ is an entourage of $Y$. Given an entourage $E$ of $X$, an $\bm{E}$\textbf{-chain} in $X$ is a finite sequence $x_1,\ldots ,x_n$ such that $(x_i,x_{i+1})\in E$ for each $i\leq n$. Inverses and concatenations of $E$-chains are defined in the obvious way. $X$ is \textbf{chain connected} if for each entourage $E$ of $X$ and any $x,y\in X$ there is an $E$-chain starting at $x$ and ending at $y$. A function $f:X\to Y$ from a uniform space $X$ has \textbf{chain lifting} if for every entourage $E$ of $X$ there is an entourage $F$ of $X$ so that for any $x\in X$, any $f(F)$-chain in $Y$ starting at $f(x)$ can be lifted to an $E$-chain in $X$ starting at $x$. The function $f$ has \textbf{uniqueness of chain lifts} if for every entourage $E$ of $X$ there is an entourage $F\subset E$ so that any two $F$-chains in $X$ starting at the same point with identical images must be equal. Showing that $f$ has unique chain lifting amounts to finding an entourage of $X$ that is \textbf{transverse} to $f$. An entourage $E_0$ is transverse to $f$ if for any $(x,y)\in E_0$ with $f(x)=f(y)$, we must have $x=y$. The function has \textbf{unique chain lifting} if it has both chain lifting and uniqueness of chain lifts. Define a function $f:X\to Y$ is a \textbf{uniform covering map} if it generates the uniform structure on $Y$ and has unique chain lifting. 

Like in the setting of paths, we wish to have homotopies of chains. Homotopies between chains were successfully defined in \cite{BP}. The following is an equivalent definition from \cite{Rips} that relies on homotopies already defined for paths. It utilizes Rips complexes which are a fundamental tool for studying chains in a uniform space. Given an entourage $E$ of $X$ the Rips complex $R(X,E)$ is the subcomplex of the full complex over $X$ whose simplices are finite $E$-bounded subsets of $X$. Any $E$-chain $x_1,\ldots ,x_n$ determines a homotopy class of paths in $R(X,E)$. Simply join successive terms $x_i,x_{i+1}$ by an edge path, i.e., a path along the edge joining $x_i$ and $x_{i+1}$. Since only homotopy classes of paths will be considered any two such paths will be equivalent. Two $E$-chains starting at the same point $x$ and ending at the same point $y$ are $\bm{E}$\textbf{-homotopic relative endpoints} if the corresponding paths in $R(X,E)$ are homotopic relative endpoints.

We wish to consider finer and finer chains in a space and therefore come to the concept of generalized paths. A \textbf{generalized path} is a collection of homotopy classes of chains $\alpha=\{[\alpha_E]\}_E$ where $E$ runs over all entourages of $X$ and for any $F\subset E$, $\alpha_F$ is $E$-homotopic relative endpoints to $\alpha_E$. Inverses and concatenations of generalized paths are defined in the obvious way. The set of generalized paths in $X$ starting at $x_0$ is denoted as $GP(X,x_0)$. We will suppress the use of a basepoint and just write $GP(X)$. $GP(X)$ is given a uniform structure generated by basic entourages defined as follows. For each entourage $E$ of $X$ let $E^*$ be the set of all pairs $(\alpha,\beta)$, $\alpha,\beta \in GP(X,x_0)$, such that $\alpha^{-1}\beta$ is $E$-homotopic to the chain $x,y$ where $x$ is the endpoint of $\alpha$ and $y$ is the endpoint of $\beta$. Call such a generalized path $\bm{E}$\textbf{-short}.

More definitions are needed in order to define a generalized uniform covering map. Suppose $f:X\to Y$ is a function between uniform spaces. This function has \textbf{approximate uniqueness of chain lifts} if for each entourage $E$ of $X$ there is an entourage $F\subset E$ such that any two $F$-chains that start at the same point and have identical images under $f$ are $E$-close. Two chains $x_1,\ldots,x_n$ and $y_1,\ldots,y_n$ are $\bm{E}$\textbf{-close} if $(x_i,y_i)\in E$ for each $i\leq n$. The function $f$ has \textbf{generalized path lifting} if for any $x\in X$, any generalized path starting at $f(x)$ lifts to a generalized path starting at $x$.

A map $f:X\to Y$ is a \textbf{generalized uniform covering map} if it generates the uniform structure on $Y$ and has chain lifting, approximate uniqueness of chain lifts, and generalized path lifting. If $f$ has complete fibers, then the requirement that it has generalized path lifting can be removed. \cite[Lemma 5.7]{Rips} In this paper we focus on generalized uniform covering maps with complete fibers since inverse limits of uniform covering maps have complete fibers (see \ref{InvLimHasCompFibers}). 

Conditions for the endpoint map $GP(X,x_0)\to X$ to be a generalized uniform covering map are introduced in \cite{Rips}. A uniform space $X$ is \textbf{uniform joinable} if any two points in $X$ can be joined by a generalized path. A uniform space $X$ is \textbf{locally uniform joinable} if for each entourage $E$ of $X$ there is an entourage $F\subset E$ such that if $(x,y)\in F$, $x$ and $y$ can be joined by an $E$-short generalized path. It is easy to see that $X$ is locally uniform joinable and chain connected if and only if $X$ is locally uniform joinable and uniform joinable. We will use the former description. The endpoint map $GP(X,x_0)\to X$ is a generalized uniform covering map if and only if $X$ is locally uniform joinable chain connected.

\section{Generalized uniform covering maps approximated by uniform covering maps}

The space $GP(X)$ was considered by Berestovskii and Plaut as the inverse limit of spaces $X_E$ for which $X_E\to X$ is realized as the projection associated with a group acting on $X_E$ that is uniformly equicontinuous and uniformly properly discontinuous. Thus, in that situation, the map $X_E\to X$ is a uniform covering map (see Section \ref{GroupActions}). Let us describe $X_E$ geometrically and see that $X_E\to X$ is a uniform covering map without appealing to group actions. 

Fix a basepoint of a uniform space $X$. Given an entourage $E$ of $X$, denote the set of $E$-homotopy classes of $E$-chains in $X$ starting at the basepoint as $X_E$. Berestovskii and Plaut considered the following uniform structure on $X_E$. For each entourage $F\subset E$, let $\widehat{F}$ be the set of all pairs $([c],[d])\in X_E\times X_E$, such that $c^{-1}d$ is $E$-homotopic relative endpoints to the edge between the endpoints of $c$ and $d$ and those endpoints are $F$-close. \footnote{Berestovskii and Plaut denote this basic entourage as $F^*$.} Then $\{\widehat{F}:F\subset E\}$ is a basis for a uniform structure on $X_E$. 

If $X$ is a chain connected uniform space then, given any entourage $E$ of $X$, the endpoint map $p_E:X_E\to X$ is a uniform covering map. It generates the uniform structure on $X$ since if $F\subset E$, $F=p_E(\widehat{F})$ (chain connectivity is used here). To see that $p_E$ has chain lifting, suppose $F\subset E$, $[c]\in X_E$, and $y\in X$ with $(p_E([c]),y)\in F$. Then $c$ concatenated with $y$ is an $E$-chain whose equivalence class is $\widehat{F}$-close to $[c]$. Finally, $\widehat{E}$ is transverse to $p_E$ so that $p_E$ has unique chain lifting.

For entourages $F\subset E$, an $F$-chain is also an $E$-chain and if two $F$-chains are $F$-homotopic relative endpoints then they are also $E$-homotopic relative endpoints. Therefore there is a map $\phi_{FE}:X_F\to X_E$ that sends an equivalence class $[c_F]_F$ in $X_F$ to the equivalence class $[c_F]_E$ in $X_E$. The corresponding inverse limit $\varprojlim X_E$ is uniformly equivalent to $GP(X)$. Therefore the endpoint map $\pi_X:GP(X)\to X$ is approximated by uniform covering maps in the sense that for any basic entourage $E^*$ of $GP(X)$, $GP(X)\to X$ factors as $GP(X)\overset{\pi_E}{\rightarrow}X_E\overset{p_E}{\rightarrow}X$ with $g_E$ a uniform covering map and $\pi_E$ having $E^*$-bounded fibers, where $\pi_E$ is the projection $([c_E])\mapsto [c_E]$. Since $GP(X)\to X$ is the model generalized uniform covering map, it makes sense to see if generalized uniform covering maps are approximated by uniform covering maps.

\begin{definition}
A map $f:X\to Y$ between uniform spaces is approximated by uniform covering maps if for every entourage $E$ of $X$, there is an entourage $F\subset E$ so that $f$ factors as $X\overset{h_F}{\to}Z_F\overset{g_F}{\to}Y$ where $g_F$ is a uniform covering map and $h_F$ has $F$-bounded fibers.
\end{definition}

Investigating when a generalized uniform covering map is approximated by uniform covering maps will be a step toward characterizing generalized uniform covering maps that are inverse limits of uniform covering maps. There is a stronger condition than approximate uniqueness of chain lifts that is necessary for a generalized uniform covering map to be approximated by uniform covering maps.

\begin{definition}
A map $f:X\to Y$ between uniform spaces has strong approximate uniqueness of chain lifts if for every entourage $E$ of $X$ there is an entourage $F\subset E$ so that any two $F$-chains starting at the same point that have identical images are $F$-close.
\end{definition}

This condition is stronger than approximate uniqueness of chain lifts because the chains are  required to be $F$-close. Whether the condition is strictly stronger is unknown.

\begin{question}
Does every generalized uniform covering map have strong approximate uniqueness of chain lifts?
\end{question}

The following shows that the condition is necessary for a map to be approximated by uniform covering maps.

\begin{proposition} \label{AppByUCMsImpliesStongAppUniq}
Suppose $f:X\to Y$ is a map between uniform spaces that is approximated by uniform covering maps. Then $f$ has strong approximate uniqueness of chain lifts.
\end{proposition}

\proof
Given an entourage $E$ of $X$, take an entourage $F\subset E$ so that $f$ factors as $X\overset{h_F}{\to}Z_F\overset{g_F}{\to}Y$ where $g_F$ is a uniform covering map and $h_F$ has $F$-bounded fibers. Let $E_0$ be an entourage of $Z_F$ that is transverse to $g_F$. Suppose there are two $h_F^{-1}(E_0)\cap F$-chains $c$ and $d$ starting at the same point with identical images under $f$. Then $h_F(c)$ and $h_F(d)$ are $E_0$-chains starting at the same point with identical images so they are identical. Since $h_F$ has $F$-bounded fibers, $c$ and $d$ are $F$-close. Of course they are also $h_F^{-1}(E_0)$-close since their images under $h_F$ are identical.
\endproof

Now we wish to mimic the factoring of $GP(X)\to X$ as $GP(X)\to X_E\to X$ for arbitrary generalized uniform covering maps with strong approximate uniqueness. Given a map $f:X\to Y$ between uniform spaces, an entourage $E$ of $X$, a set $A\subset X$, and an $x\in A$, the $\bm{E}$\textbf{-component} of $x$ in $A$ is the set of all $y\in A$ that can be joined to $x$ by an $E$-chain in $A$. Given a map $f:X\to Y$ and an entourage $E$ of $X$, let $X_f/E$ be the set obtained by identifying the $E$-components of the fibers of $f$ to a point. Let $q_E:X\to X_f/E$ denote the quotient function associated with the identification. Note that if $F$ is an entourage such that two $F$-chains with the same images under $f$ are $F$-close then $q_F(x)=q_F(y)$ if and only if $f(x)=f(y)$ and $(x,y)\in F$. Define a map $g_F:X_f/F\to Y$ to send $q_F(x)$ to $f(x)$. Note it is well defined since $q_F(x)=q_F(x')$ implies $f(x)=f(x')$.

To give $X_f/E$ a uniform structure, of course we could just pull back the uniform structure of $Y$, but then $g_E$ will not be a uniform covering map. In fact we will only consider a uniform structure on $X_f/F$ for basic entourages $F$ of $X$.  

\begin{lemma} \label{q_FGeneratesStructure}
Suppose a map $f:X\to Y$ has chain lifting and $F$ is an entourage of $X$ such that two $F$-chains that start at the same point and have the same images under $f$ are $F$-close. Then $q_F$ generates a uniform structure on $X_f/F$.
\end{lemma}

\proof
It suffices to show that $q_F$ has chain lifting. Given an entourage $D$ of $X$ choose $H\subset D\cap F$ so that $f(H)$-chains in $Y$ lift to $D\cap F$-chains in $X$. Suppose $x,y\in X$ with $(q_F(x),q_F(y))\in q_F(H)$. Then $q_F(x)=q_F(x')$ and $q_F(y)=q_F(y')$ for some $(x',y')\in H$. Then $(f(x),f(y))\in f(H)$ so there is a $y''\in X$ with $f(y'')=f(y)$ and $(x,y'')\in D\cap F$. But then $x,x,y''$ and $x,x',y'$ are two $F$-chains with identical images under $f$ so $(y',y'')\in F$. Then $q_F(y'')=q_F(y')=q_F(y)$.  
\endproof

\begin{proposition} \label{fAppByUCM}
Suppose $f:X\to Y$ generates the uniform structure on $Y$, has chain lifting, and has strong approximate uniqueness of chain lifts. Then $f$ is approximated by uniform covering maps.
\end{proposition}

\proof
By \ref{q_FGeneratesStructure}, for each entourage $E$ of $X$ there is an entourage $F\subset E$ so that $f$ factors as $X\overset{q_F}{\to}X_f/F\overset{g_F}{\to}Y$ where $q_F$ has $F$-bounded fibers. Since $f$ has chain lifting and $g_F$ generates the uniform structure on $X_f/F$, $g_F$ has chain lifting. \cite{Subgroups} Note that $g_F$ generates the uniform structure on $Y$ since $g_F(q_F(G))=f(G)$. It remains to show that there is an entourage of $X_f/F$ that is transverse to $g_F$. Suppose $(q_F(x),q_F(y))\in q_F(F)$ with $f(x)=f(y)$. Then $q_F(x)=q_F(x')$ and $q_F(y)=q_F(y')$ for some $(x',y')\in F$. Then $x,x',y',y$ is an $F$-chain with $f(x)=f(x')=f(y')=f(y)$ so $q_F(x)=q_F(y)$.
\endproof

As justification for the definition of $X_f/E$, we will show that if $X$ is locally uniform joinable chain connected then $GP(X)_{\pi_X}/E^*$ is uniformly equivalent to $\pi_E(GP(X))$. First notice that there is a uniform structure on $GP(X)_{\pi_X}/E^*$ for any basic entourage $E^*$ of $GP(X)$.

\begin{lemma} \label{GP(X)AppUniq}
For any basic entourage $E^*$ of $GP(X)$, any two $E^*$-chains that start at the same point and have identical images are $E^*$-close.
\end{lemma}

\proof
Suppose $\alpha_1,\ldots,\alpha_n$ and $\beta_1,\ldots,\beta_n$ are $E^*$-chains with identical images under $\pi$ and $\alpha_1=\beta_1$. It suffices to show that if $(\alpha_i,\beta_i)\in E^*$ then $(\alpha_{i+1},\beta_{i+1})\in E^*$. If $(\alpha_i,\beta_i)\in E^*$ then $[\alpha_{i+1}^{-1} \beta_{i+1}]_E=[\alpha_{i+1}^{-1} \alpha_i \beta_i^{-1} \beta_{i+1}]_E$ which is $E$-homotopic to the constant chain at the endpoint of $\alpha_{i+1}$.
\endproof

\begin{proposition}
If $X$ is locally uniform joinable chain connected then \newline $GP(X)_{\pi_X}/E^*$ is uniformly equivalent to $\pi_E(GP(X))$.
\end{proposition}

\proof
First note there is a bijective correspondence between the two sets. Let $i:GP(X)_{\pi_X}/E^*\to X_E$ send $q_{E^*}(\alpha)$ to $\alpha_E$. The function is well defined and injective. Of course it need not be surjective since $X_E$ contains all $E$-chains starting at $x_0$ while $i(GP(X)_{\pi_X}/E^*)$ only contains $E$-chains that are terms of generalized paths. But it therefore does map onto $\pi_E(GP(X))$. To see that $i$ generates the uniform structure on $\pi_E(GP(X))$, first note that $i$ is uniformly continuous since for an entourage $F\subset E$, $q_{F^*}\subset i^{-1}(\widehat F)$. Now consider a basic entourage $q_{E^*}(F^*)$ of $GP(X)_{\pi_X}/E^*$. Let $H\subset E$ be an entourage of $X$ so that if $(\pi_X(\alpha),y)\in H$ there is a $\beta\in GP(X)$ with $(\alpha,\beta)\in E^*\cap F^*$ and $\pi_X(\beta)=y$. Suppose $\alpha,\beta\in GP(X)$ with $(\pi_E(\alpha),\pi_E(\beta))\in \widehat{H}$. Set $\pi_X(\alpha)=x$ and $\pi_X(\beta)=y$. Then $(x,y)\in H$ so there is a $\beta'\in GP(X)$ with $\pi_X(\beta')=y$ and $(\alpha,\beta')\in E^*\cap F^*$. Then $\alpha, \beta$ and $\alpha, \beta'$ are two $E^*$ chains with identical images so $(\beta, \beta ')\in E^*$. Therefore $q_{E^*}(\beta)=q_{E^*}(\beta ')$ so $(q_{E^*}(\alpha),q_{E^*}(\beta))\in q_{E^*}(F^*)$. Notice $i(q_{E^*}(\alpha),q_{E^*}(\beta))=(\pi_E(\alpha),\pi_E(\beta))$.
\endproof

\section{Inverse limits of generalized uniform covering maps}

Now inverse limits are investigated. Recall that $GP(X)\to X$ is the inverse limit of the maps $X_E\to X$ and these maps are uniform covering maps provided $X$ is chain connected. Now $GP(X)\to X$ is only a generalized covering map if $X$ is also uniform joinable, which is equivalent to the inverse system $\{X_E,\phi_{FE}\}$ being strong Mittag-Leffler. An inverse system $\{X_\alpha,\phi_{\beta \alpha}\}$ satisfies the strong Mittag-Leffler condition if for each $\alpha$ there is a $\beta<\alpha$ such that $\phi_{\beta \alpha}(X_\beta)\subset \pi_\alpha(\varprojlim X_\alpha)$ \cite{Rips}. Note that this inclusion implies that the sets are in fact equal. Therefore we consider inverse limits of strong Mittag-Leffler inverse systems of uniform covering maps.

The strong Mittag-Leffler condition is assumed for many of the results of this section but it is not needed for strong approximate uniqueness of chain lifts to be preserved by inverse limits.

\begin{lemma} \label{InvLimStrongAppUniq}
Suppose there are inverse systems of uniform spaces $\{X_\alpha,\phi_{\beta \alpha}\}$ and \linebreak $\{Y_\alpha,\psi_{\beta \alpha}\}$ with compatible maps $f_\alpha:X_\alpha \rightarrow Y_\alpha$. Set $f=\ilim f_\alpha$.

\begin{itemize}
\item[1.] If each $f_\alpha$ has approximate uniqueness of chain lifts then $f$ has approximate uniqueness of chain lifts.
\item[2.] If each $f_\alpha$ has strong approximate uniqueness of chain lifts then $f$ has strong approximate uniqueness of chain lifts.
\end{itemize}
\end{lemma}

\proof
1. Consider a basic entourage $\pi_\alpha ^{-1}(E_\alpha)$ of $\varprojlim (X_\alpha)$. Now there is an entourage $F_\alpha \subset E_\alpha$ so that any two $F_\alpha$-chains in $X_\alpha$ starting at the same point that have identical images under $f_\alpha$ are $E_\alpha$-close. Suppose $c$ and $d$ are two $\pi_\alpha^{-1}(F_\alpha)$-chains in $\varprojlim (X_\alpha)$ starting at the same point that have identical images under $f$. Then $\pi_\alpha(c)$ and $\pi_\alpha(d)$ are two $F_\alpha$-chains in $X_\alpha$ starting at the same point with identical images under $f_\alpha$ so they are $E_\alpha$-close. Therefore $c$ and $d$ are $\pi_\alpha ^{-1}(E_\alpha)$-close.

2. The proof is similar to the proof of 1.
\endproof

The following shows that having complete fibers is a necessary condition for a generalized uniform covering map to be the inverse limit of uniform covering maps.

\begin{lemma} \label{InvLimHasCompFibers}
Suppose there are inverse systems $\{X_\alpha,\phi_{\beta \alpha}\}$ and $\{Y_\alpha,\psi_{\beta \alpha}\}$ of uniform spaces and compatible maps $f_\alpha:X_\alpha \rightarrow Y_\alpha$ that have complete fibers. Suppose each $X_\alpha$ is Hausdorff. Then $f=\varprojlim f_\alpha$ has complete fibers.
\end{lemma}

\proof
Suppose $(y_\alpha)\in \varprojlim Y_\alpha$. Notice $f^{-1}((y_\alpha))$ is identical to $\varprojlim f_\alpha^{-1}(y_\alpha)$ as sets and uniform spaces. Since the inverse limit of complete Hausdorff spaces is complete \cite{B}, $f^{-1}((y_\alpha))$ is complete.
\endproof

Now consider a strong Mittag-Leffler inverse system.

\begin{lemma} \label{InverseLimitsOfUCMs}
Suppose there is a strong Mittag-Leffler inverse system of uniform spaces $\{X_\alpha,\phi_{\beta \alpha}\}$ with compatible maps $f_\alpha:X_\alpha \rightarrow
Y$. Set $f=\underset{\leftarrow}{\lim }f_\alpha$.
\begin{itemize}
\item[1.] If each $f_\alpha$ generates the uniform structure on $Y$ then $f$
generates the uniform structure on $Y$.
\item[2.] If each $f_\alpha$ has chain lifting then $f$ has chain lifting.
\end{itemize}
\end{lemma}

\proof
Consider a basic entourage $\pi_\alpha^{-1}(E_\alpha)$ of $\varprojlim X_\alpha$. Let $\beta <\alpha$ so that $\phi_{\beta \alpha}(X_\beta)=\pi_\alpha(\varprojlim X_\alpha)$.

1. First note that $f$ is uniformly continuous since each $f_\alpha$ is. To see that $f_\beta \phi_{\beta \alpha}^{-1}(E_\alpha)\subset
f(\pi_\alpha^{-1}(E_\alpha))$, suppose $(x,y)\in f_\beta \phi_{\beta \alpha}^{-1}(E_\alpha)$. Then $(x,y)=f_\beta (x_\beta^\prime,y_\beta^\prime)$ for some $(x_\beta^\prime,y_\beta^\prime)\in \phi_{\beta \alpha }^{-1}(E_\alpha)$. Now $\phi_{\beta \alpha}(x_\beta^\prime,y_\beta^\prime)=\pi_\alpha((x_\alpha),(y_\alpha))$ for some $(x_\alpha),(y_\alpha)\in \varprojlim X_\alpha$. Note \newline $f((x_\alpha),(y_\alpha))=(x,y)$. Also $\pi_\alpha((x_\alpha),(y_\alpha))=\phi_{\beta
\alpha}(x_\beta^\prime,y_\beta^\prime)\in E_\alpha$ so \newline $(x,y)\in f(\pi_\alpha^{-1}(E_\alpha))$.

2. There is an entourage $F_\beta$ of $X_\beta$ so that an if $(f_\beta(x),y)\in f_\beta(F_\beta)$, there is a $y'\in X_\beta$ with $f_\beta(y')=y$ and $(x,y')\in \phi_{\beta \alpha}^{-1}(E_\alpha)$. Suppose $(f((x_\alpha)),y)\in f(\pi_\beta^{-1}(F_\beta))$. Then $(f_\beta(x_\beta),y)\in f_\beta(F_\beta)$ so there is a $y'_\beta\in X_\beta$ with $f_\beta(y'_\beta)=y$ and $(x_\beta,y'_\beta)\in \phi_{\beta \alpha}^{-1}(E_\alpha)$. Now $\phi_{\beta \alpha}(y'_\beta)=\pi_\alpha((y_\alpha))$ for some $(y_\alpha)\in \varprojlim X_\alpha$. Note $f((y_\alpha))=y$ and $((x_\alpha),(y_\alpha))\in \pi_\alpha^{-1}(E_\alpha)$.
\endproof

Notice it is unclear how to generalize this result to an inverse system $\{Y_\alpha,\psi_{\beta \alpha}\}$ and maps $f_\alpha:X_\alpha \rightarrow Y_\alpha$. There is even a problem showing that the resulting map $f$ is surjective.

\begin{proposition}
Suppose there is a strong Mittag-Leffler inverse system $\{X_\alpha ,\phi_{\beta \alpha}\}$ of Hausdorff uniform spaces and compatible generalized uniform covering maps $f_\alpha:X_\alpha \to Y$ with complete fibers. Then the inverse limit $f=\varprojlim f_\alpha$ is a generalized uniform covering map.
\end{proposition}

\proof
By \ref{InvLimHasCompFibers} the inverse limit has complete fibers. By \ref{InverseLimitsOfUCMs} it generates the uniform structure on $Y$ and has chain lifting. By \ref{InvLimStrongAppUniq} it has approximate uniqueness of chain lifts.
\endproof

In particular, the inverse limit of a strong Mittag-Leffler inverse system of uniform covering maps over a Hausdorff space is a generalized uniform covering map.

Now we wish the express a generalized uniform covering map as the inverse limit of uniform covering maps. According to \ref{AppByUCMsImpliesStongAppUniq}, strong approximate uniqueness is a necessary condition. Then, according to \ref{fAppByUCM}, given a generalized uniform covering map $f:X\to Y$ that has strong approximate uniqeness of chain lifts, $X$ has a basis of entourages such that for each basic entourage $E$, $q_E:X\to X_f/E$ has $E$-bounded fibers and $g_E:X_f/E\to Y$ is a uniform covering map. Consider the inverse system $\{X_f/E,\phi_{FE}\}$ given by this basis where for basic entourages $F\subset E$ of $X$, $\phi_{FE}:X_f/F\to X_f/E$ is defined to send $q_F(x)$ to $q_E(x)$. These functions are well-defined since each equivalence class of $X_f/F$ is contained in a single equivalence class of $X_f/E$. The functions are also uniformly continuous and compatible with the maps $\{q_E\}$ and $\{g_E\}$.

\begin{proposition} \label{X=limX_f/E}
Suppose $f:X\to Y$ is a generalized uniform covering map that has strong approximate uniqueness. Consider the inverse system $\{X_f/E,\phi_{FE}\}$ given by the basis as above and set $q=\varprojlim q_E$. Suppose $X$ is Hausdorff. Then $q$ is a uniform embedding.  If $f$ has complete fibers then $q$ is a uniform equivalence.
\end{proposition}

\proof
Since each $q_E$ is uniformly continuous $q$ is as well. To see that $q$ is injective, suppose $x,y\in X$ with $q(x)=q(y)$. Since $q_E(x)=q_E(y)$ for all basic of entourages $E$ of $X$, $(x,y)\in E$ for each basic entourage $E$. Then since $X$ is Hausdorff, $x=y$. To see that $q$ is a uniform embedding, suppose $E$ is a basic entourage of $X$. Choose a basic entourage $F$ of $X$ so that $F^3\subset E$. Then to show that $q$ is a uniform embedding it suffices to show that for any $x,y\in X$ with $(q(x),q(y))\in \pi_F^{-1}(q_F(F))$, $x,y\in E$ where $\pi_F$ is the projection from $\varprojlim X_f/E$ to $X_f/F$. Suppose $(q(x),q(y))\in \pi_F^{-1}(q_F(F))$. Then there are $x',y'\in X$ with $(x,x'),(x',y'),(y',y)\in F$. Therefore $(x,y)\in E$.

Finally, to see that $q$ is surjective and therefore a uniform equivalence if $f$ has complete fibers, suppose $(q_E(x_E))\in \varprojlim X_f/E$. Set $y=f(x_E)$. Note $y$ is independent of the entourage $E$ and the choice of representative $x_E$. For each basic entourage $E$ let $A_E\subset X$ be the equivalence class of $x_E$. Then $\{A_E\}$ is a Cauchy filter base in the fiber $f^{-1}(y)$ so it has a limit point $x$. Now, given an entourage $E$ of $X$ there is an entourage $F$ so that $A_F\subset B(x,E)$. Then $A_{F\cap E}\subset B(x,E)$ so $q_E(x_{F\cap E})=q_E(x)$. But $q_E(x_{F\cap E})=q_E(x_E)$ so $q(x)=(q_E(x_E))$.
\endproof

\begin{theorem}
Suppose $f:X\to Y$ is a map between uniform spaces with $Y$ Hausdorff. Then the following are equivalent.
\begin{itemize}
\item[1.] There is a strong Mittag-Leffler inverse system $\{X_\alpha,\phi_{\beta \alpha}\}$ and uniform covering maps $f_\alpha:X_\alpha \to Y$ with $f=\varprojlim f_\alpha$.
\item[2.] $X$ is Hausdorff and $f$ generates the uniform structure on $Y$, has chain lifting, has strong approximate uniqueness of chain lifts, and has complete fibers.
\end{itemize}
\end{theorem}

\proof
1. $\implies$ 2. By \ref{InverseLimitsOfUCMs}, $f$ generates the uniform structure on $Y$ and has chain lifting. By \ref{InvLimStrongAppUniq} or \ref{AppByUCMsImpliesStongAppUniq} $f$ has strong approximate uniqueness. Notice that each $X_\alpha$ is Hausdorff as $Y$ is Hausdorff and $f_\alpha$ is a uniform covering map. Then $f$ has complete fibers by \ref{InvLimHasCompFibers}. Notice that the inverse limit of Hausdorff spaces is Hausdorff.

2. $\implies$ 1. By \ref{X=limX_f/E}, $f$ is uniformly equivalent to the map $\varprojlim X_f/E\to Y$. Notice the inverse system is strong Mittag-Leffler since $X\to X_f/E$ is surjective for each basic entourage $E$.
\endproof

\section{Inverse limits of regular generalized uniform covering maps} \label{GroupActions}

Since generalized uniform covering maps can be characterized as the inverse limit of uniform covering maps (provided that the map has strong approximate uniqueness of chain lifts), we now attempt to characterize generalized uniform covering maps induced by group actions as inverse limits of uniform covering maps induced by group actions. 

``Group Actions and Covering Maps in the Uniform Category'' \cite{GroupActions} considered group actions that induce uniform covering maps and generalized uniform covering maps. Unless otherwise noted, the following definitions and results are from that paper. 

Suppose a group $G$ acts on a uniform space $X$. The action is \textbf{neutral} \cite{J} if for each entourage $E$ of $X$ there is an entourage $F$ of $X$ such that if $(x,gy)\in F$ there is an $h\in G$ with $(hx,y)\in E$. Note that if $G$ acts neutrally on $X$, then the projection $p:X\to X/G$ has chain lifting and generates a uniform structure on $X/G$. This is the structure that will always be considered. The action is \textbf{uniformly properly discontinuous} if there is an entourage $E_0$ of $X$ such that if $(x,gx)\in E_0$ for some $x\in X$ then $g=1$. If $G$ is acts faithfully on a chain connected uniform space $X$ then the action is small scale uniformly continuous and the projection $p:X\to X/G$ is a uniform covering map if and only if the action is neutral and uniformly properly discontinuous.

Given an entourage $E$ of $X$ define $\bm{G_E}$ to be the subgroup generated by $\{g\in G:(x,gx)\in E \text{ for some } x\in X\}$ \cite{BP}. The action is \textbf{small scale uniformly continuous} if for each entourage $E$ of $X$ there is an entourage $F$ of $X$ so that for each $g\in G_F$, $g^{-1}(E)$ is an entourage of $X$. The action is \textbf{small scale uniformly equicontinuous} if for each entourage $E$ of $X$ there is an entourage $F$ of $X$ such that for each $g\in G_F$, $F\subset g^{-1}(E)$. The action is \textbf{uniformly equicontinuous} \cite{J} if for each entourage $E$ of $X$ there is an entourage $F$ of $X$ such that for each $g\in G$, $F\subset g^{-1}(E)$. Equivalently, $X$ has a basis of $G$-invariant entourages. An entourage $E$ is $\bf{G}$\textbf{-invariant} if for each $g\in G$, $gE=E$.  Finally, the action has \textbf{small scale bounded orbits} if for each entourage $E$ of $X$ there is an entourage $F$ of $X$ such that the action of $G_F$ on $X$ has $E$-bounded orbits. If $X$ is chain connected and the projection $p:X\to X/G$ has complete fibers then the action is small scale uniformly equicontinuous and $p$ is a generalized uniform covering map if and only if the action is neutral and has small scale bounded orbits.

We now consider inverse limit of neutral and uniformly properly discontinuous group actions.

\begin{lemma} \label{InvLimOfNeutral}
Suppose there is an inverse system of groups $\{G_\alpha,\psi_{\beta \alpha} \}$, an inverse system of uniform spaces $\{X_\alpha,\phi_{\beta \alpha} \}$, and compatible actions of $G_\alpha$ on $X_\alpha$ that are neutral. Suppose the inverse system $\{G_\alpha,\psi_{\beta \alpha} \}$ is strong Mittag-Leffler. Set $G=\varprojlim G_\alpha$ and $X=\varprojlim X_\alpha$. Then the induced action of $G$ on $X$ is neutral.
\end{lemma}

\proof
Let $\pi_\alpha^{-1}(E_\alpha)$ be a basic entourage of $X$. There is a $\beta<\alpha$ with $\psi_{\beta \alpha}(G_\beta)=\pi_\alpha (G)$. Now there is an entourage $F_\beta \subset \phi_{\beta \alpha}^{-1}(E_\alpha)$ so that if $(x,hy)\in F_\beta$ then there is a $g\in G_\beta$ with $(gx,y)\in \phi_{\beta \alpha}^{-1}(E_\alpha)$. Suppose $((x_\alpha),(h_\alpha)(y_\alpha))\in \pi_\beta^{-1}(F_\beta)$. Then $(x_\beta,h_\beta y_\beta)\in F_\beta$ so there is a $g'_\beta \in G_\beta$ with $(g'_\beta x_\beta,y_\beta)\in \phi_{\beta \alpha}^{-1}(E_\alpha)$. Now there is a $(g_\alpha)\in G$ with $g_\alpha=\psi_{\beta \alpha}(g'_\beta)$. Then $((g_\alpha)(x_\alpha),(y_\alpha))\in \pi_\alpha^{-1}(E_\alpha)$.
\endproof

\begin{lemma} \label{InvLimOfSmScBddOrbits}
Suppose there is an inverse system of groups $\{G_\alpha,\psi_{\beta \alpha} \}$, an inverse system of uniform spaces $\{X_\alpha,\phi_{\beta \alpha} \}$, and compatible actions of $G_\alpha$ on $X_\alpha$ that have small scale bounded orbits. Set $G=\varprojlim G_\alpha$ and $X=\varprojlim X_\alpha$. Then the induced action of $G$ on $X$ has small scale bounded orbits.
\end{lemma}

\proof
Let $\pi_\alpha^{-1}(E_\alpha)$ be a basic entourage of $X$. Let $F_\alpha$ be an entourage of $X_\alpha$ so that orbits of the action of $G_{\alpha \:F_\alpha}$ on $X_\alpha$ are $E_\alpha$-bounded. Suppose $(g_\alpha)\in G_{\pi_\alpha^{-1}(F_\alpha)}$ and $(x_\alpha)\in X$. Then $g_\alpha \in G_{\alpha \:F_\alpha}$ so $(x_\alpha,g_\alpha x_\alpha)\in E_\alpha$. Therefore $((x_\alpha),(g_\alpha)(x_\alpha))\in \pi_\alpha^{-1}(E_\alpha)$.
\endproof

In particular the inverse limit of uniformly properly discontinuous actions has small scale bounded orbits.

Given a strong Mittag Leffler inverse system of groups $\{G_\alpha,\psi_{\beta \alpha} \}$, an inverse system of Hausdorff uniform spaces $\{X_\alpha,\phi_{\beta \alpha} \}$, and compatible actions of $G_\alpha$ on $X_\alpha$ that are neutral and properly discontinuous, we wish to have the projection $p:X\to X/G$ be a generalized uniform covering map where $G=\varprojlim G_\alpha$ and $X=\varprojlim X_\alpha$. By \ref{InvLimOfNeutral} the action is neutral and by \ref{InvLimOfSmScBddOrbits} the action has small scale bounded orbits. It remains to have that the projection $X\to X/G$ has complete fibers. For this proof we need a strong assumption on the inverse system $\{G_\alpha,\psi{\beta \alpha} \}$. In particular we need it to be a countable sequence with $\ilim^1 G_i=1$. 

Given an inverse sequence of groups $\{G_i,\psi_i\}$, consider the function $f:\Pi G_i\to \Pi G_i$ that sends $(g_i)$ to $(g_i\psi_{i+1}(g_{i+1})^{-1})$. Then $\ilim G_i$ and $\ilim^1 G_i$ can be defined as sets for which $1\to \ilim G_i\hookrightarrow \Pi G_i\stackrel{f}{\to} \Pi G_i\hookrightarrow\ilim^1 G_i\to 1$ is an exact sequence. Of importance to us is the case when $\ilim^1 G_i=1$. Notice $\ilim^1 G_i=1$ if and only if $f$ is surjective. Therefore we have the following.

\begin{lemma}
Given an inverse sequence of groups $\{G_i,\psi_i\}$, $\ilim^1 G_i=1$ if and only if for each sequence $(g_i)\in \Pi G_i$, there is a sequence $(h_i)\in \Pi G_i$ with $g_i=h_i\psi_{i+1}(h_{i+1})^{-1}$ for all $i\geq 1$.
\end{lemma}

Notice if each $\psi_i$ is surjective then $\ilim^1 G_i=1$. Indeed, given $(g_i)\in \Pi G_i$, set $h_1=1$. Then for each $i\geq 1$ choose $h_{i+1}$ so that $\psi_{i+1}(h_{i+1})=g_i^{-1}h_i$. 

It is helpful to also have the following.

\begin{lemma} \label{Lim1}
If $\ilim^1 G_i=1$ then for each sequence $(g_i)\in \Pi G_i$, there is a sequence $(h_i)\in \Pi G_i$ with $g_i=\psi_{i+1}(h_{i+1})^{-1}h_i$ for all $i\geq 1$.
\end{lemma}

\proof
Given the sequence $(g_i)$, there is a sequence $(k_i)$ with $g_i^{-1}=k_i\psi_{i+1}^{-1}(k_{i+1})$ for each $i$. Then $g_i=\psi_{i+1}(k_{i+1})k_i^{-1}$. Set $h_i=k_i^{-1}$. Then $g_i=\psi_{i+1}(h_{i+1}^{-1})h_i=\psi_{i+1}(h_{i+1})^{-1}h_i$.
\endproof

\begin{proposition} \label{InverseLimitGroupActions}
Suppose there is an inverse sequence of uniform spaces $\{X_i,\phi_i\}$ and a Mittag-Leffler inverse sequence of groups $\{G_i,\psi_i\}$ with compatible neutral and free actions of $G_i$ on $X_i$. Let $X=\ilim X_i$ and $G=\ilim G_i$. Suppose that $\ilim^1 G_i=1$. Then $\ilim (X_i/G_i)$ is uniformly equivalent to $X/G$.
\end{proposition}

\proof
We assume the actions are neutral so that there are uniform structures on $X_i/G_i$. Notice the action of $G$ on $X$ is neutral by \ref{InvLimOfNeutral} so that there is a uniform structure on $X/G$. Given $[(x_i)]\in X/G$, define $f:X/G\to\ilim (X_i/G_i)$ to send $[(x_i)]$ to $([x_i])$. Notice it is well defined. To see that it is injective, suppose $[(x_i)],[(y_i)]\in X/G$ with $([x_i])=([y_i])$. Then for each $n$ there is a $g_n\in G_n$ so that $x_n=g_n y_n$. It suffices to show that $(g_i)$ is a thread. But $\phi_n(x_n)=\phi_n(g_n y_n)=\psi_n(g_n) \phi_n(y_n)=\psi_n(g_n) y_{n-1}$ and $\phi_n(x_n)=x_{n-1}=g_{n-1} y_{n-1}$. Then since the action is free we have $\psi_n(g_n)=g_{n-1}$.

Now we will see that $f$ is uniformly continuous. For each $n$, let $q_n:X_n\to X_n/G_n$ be the quotient map associated with the action of $G_n$ on $X_n$, $\pi_n:\ilim X_i\to X_n$ be the projection to $X_n$, and $\Pi_n:\ilim (X_i/G_i)\to X_n/G_n$ be the projection to $X_n/G_n$. Let $Q:X\to X/G$ be the quotient map associated with the action of $G$ on $X$. Given an entourage $E_n$ of $X_n$, let us see that $f(Q(\pi_n^{-1}(E_n)))\subset\Pi_n^{-1}(q_n(E_n))$. Suppose $([(x_i)],[(y_i)])\in Q(\pi_n^{-1}(E_n))$. Then $(x_i)=(g_i) (x'_i)$ and $(y_i)=(h_i) (y'_i)$ with $(x'_n,y'_n)\in E_n$ for some $(g_i),(h_i)\in\ilim G_i$. Then $(([x_i]),([y_i]))\in\Pi_n^{-1}(q_n(E_n))$.

Finally, let us see that set $f(Q(\pi_n^{-1}(E_n)))$ is an entourage of $\ilim (X_i/G_i)$ so that $f$ is a uniform equivalence. Choose $m>n$ so that if $g_m\in G_m$, there is an $(h_i)\in\ilim G_i$ with $h_n=\psi_{mn}(g_m)$ where $\psi_{mn}=\psi_{n+1}\circ\psi_{n+2}\circ\cdots\circ\psi_{m-1}\circ\psi_m$. Suppose $(([x_i]),([y_i]))\in\Pi_n^{-1}(\phi_{mn}^{-1}(E_n))$. Without loss of generality we can assume $(x_i,y_i)\in\phi_{mn}^{-1}(E_n)$. We will define $(x'_i),(y'_i)\in\ilim X_i$ with $([(x'_i)],[(y'_i)])\in Q(\pi_n^{-1}(E_n))$ and $([x'_i])=([x_i])$ and $([y'_i])=([y_i])$. 

Since $([x_i])$ is a thread, for each $i$ there is a $g_i\in G_i$ with $\phi_{i+1}(x_{i+1})=g_i x_i$. Because $\ilim^1 G_i=1$ there is a sequence $\{k_i\}$ with $g_i=\psi_{i+1}(k_{i+1})^{-1} k_i$ for each $i$ (see \ref{Lim1}). Define $x'_i=k_i x_i$. Then $\phi_{i+1}(x'_{i+1})=\psi_{i+1}(k_{i+1})\phi_{i+1}(x_i)=k_i g_i^{-1}g_i x_i=x'_i$ so $(x'_i)\in\ilim X_i$. Notice $([x'_i])=([x_i])$. Define $(y'_i)\in\ilim X_i$ analogously.

We need $([(x'_i)],[(y'_i)])\in Q(\pi_n^{-1}(E_n))$. Now $(\phi_{mn}(x_m),\phi_{mn}(y_m))\in E_n$. Let us calculate $\phi_{mn}(x_m)$. We have \newline $\phi_{mn}(x_m)= \psi_{m-1\ n}(g_{m-1})\psi_{m-2\ n}(g_{m-2}) \cdots \psi_{n+2\ n}(g_{n+2})\psi_{n+1}g_{n+1} g_n x_n$. Now $g_i=\psi_{i+1}(k_{i+1})^{-1} k_i$ so 
$\psi_{in}(g_i)=\psi_{in}(\psi_{i+1}(k_{i+1})^{-1} k_i)= 
\psi_{i+1\ n}(k_{i+1})^{-1}\psi_{in}(k_i)$. Therefore $\psi_{m-1\ n}(g_{m-1})\psi_{m-2\ n}(g_{m-2})\cdots\psi_{n+1\ n}(g_{n+1})$ is a telescoping product and is equal to $\psi_{mn}(k_m)^{-1} k_n x_n$. Now there is $(h_i)\in\ilim G_i$ with $h_n=\psi_{mn}(k_m)^{-1}$ so $\pi_n((h_i)(x'_i))=\phi_{mn}(x_m)$. Similarly there is $(l_i)\in\ilim G_i$ with $\pi_n((l_i)(y'_i))=\phi_{mn}(y_m)$.
\endproof

\begin{proposition} \label{InverseLimitOfRegUCMs}
Suppose there is a Mittag-Leffler inverse sequence of groups $\{G_i,\psi_i\}$, an inverse sequence of Hausdorff uniform spaces $\{X_i,\phi_i\}$, and compatible actions of $G_i$ on $X_i$ that are neutral and have small scale bounded orbits. Suppose the projections associated with the actions of $G_i$ on $X_i$ have complete fibers. Suppose $\ilim^1 G_i=1$. Set $G=\ilim G_i$ and $X=\ilim X_i$. Then the projection $p:X\to X/G$ associated with the induced action of $G$ on $X$ is a generalized uniform covering map.
\end{proposition}

\proof
By \ref{InvLimOfNeutral} the induced action is neutral. By \ref{InvLimOfSmScBddOrbits} the induced action has small scale bounded orbits. Therefore it is enough to show that $p$ has complete fibers. Notice the action is free by \cite[Corollary 3.19]{GroupActions}. Therefore \ref{InverseLimitGroupActions} applies so the inverse limit uniform structure on $X/G$ coincides with the structure generated by $p$. But then $p$ has complete fibers by \ref{InvLimHasCompFibers}.
\endproof

In particular the inverse limit of neutral and uniformly properly discontinuous actions induces a generalized uniform covering map provided the inverse sequence of groups is Mittag-Leffler and has trivial $\ilim^1$.

Now, given an action of a group $G$ on a uniform space $X$ that induces a generalized uniform covering map $p:X\to X/G$ we wish to express that action as the inverse limit of neutral and properly discontinuous actions.

\begin{lemma} \label{BasisOfPropDiscActions}
Suppose a group $G$ acts faithfully and uniformly equicontinuously on a uniform space $X$. Then for each entourage $E$ of $X$ there is an entourage $F\subset E$ such that the induced action of $G/G_F$ on $X/G_F$ is uniformly properly discontinuous.
\end{lemma}

\proof
First note that an action is uniformly equicontinuous if and only if there is a basis of invariant entourages \cite{GroupActions}. An entourage $F$ of $X$ is $G$-invariant if $(gx,gy)\in F$ for each $(x,y)\in F$ and $g\in G$. Given an entourage $E$ of $X$, choose $F\subset E$ to be $G$-invariant. Notice since $F$ is invariant, $G_F$ is normal in $G$ \cite{P}. Then $G/G_F$ is a group and its action on $X/G_F$ is well defined. Let $p_F:X\to X/G_F$ be the projection associated with the action of $G_F$ on $X$. Since the action of $G$ on $X$ is uniformly equicontinuous, so is the action of $G_F$ on $X$ and that action is therefore neutral. Therefore $p_F$ generates a uniform structure on $X/G_F$. Suppose $([x],[g][x])\in p_F(F)$ for some $[x]\in X/G_F$ and $[g]\in G/G_F$. Then $(g_Fx,h_Fgx)\in H$ for some $g_F,h_F\in G_F$. Then $(x,g_F^{-1}h_Fgx)\in F$ so $g_F^{-1}h_Fg\in G_F$ and $g\in G_F$. Therefore $[x]=[g][x]$. Notice the action of $G/G_F$ on $X/G_F$ is faithful since $G_F$ is precicely the stabilizer of the action of $G$ on $X/G_F$. Therefore $[g]=e_{G/G_F}$.
\endproof

Notice the action of $G/G_F$ on $X/G_F$ in the above lemma is also uniformly equicontinuous. Therefore if a group $G$ acts faithfully and uniformly equicontinuously on a uniform space $X$, $X$ has a basis of entourages so that for each basic entourage $E$, the action of $G/G_E$ on $X/G_E$ induces a uniform covering map.  Consider the inverse systems $\{G/G_E,\psi_{FE}\}$ and $\{X/G_E,\phi_{FE}\}$ where for $F\subset E$, $\psi_{FE}$ sends the equivalence class of $g$ in $G/G_F$ to the equivalence class of $g$ in $G/G_E$ and $\phi_{FE}$ sends the equivalence class of $x$ in $X/G_F$ to the equivalence class of $x$ in $X/G_E$. The functions are well defined, $\psi_{FE}$ is a group homomorphism, $\phi_{FE}$ is uniformly continuous, and these maps are compatible with the actions of $G/G_E$ on $X/G_E$.

\begin{proposition} \label{G=limG/G_E}
Suppose a group $G$ acts faithfully and uniformly equicontinuously on a uniform space $X$. Suppose $X$ is Hausdorff and the action has small scale bounded orbits. Consider the inverse system $\{G/G_E,\psi_{FE}\}$ given by the basis as above. Then $G$ is isomorphic to a subgroup of $\varprojlim G/G_E$. If the projection $X\to X/G$ induced from the action of $G$ on $X$ has complete fibers then $G$ is isomorphic to $\varprojlim G/G_F$.
\end{proposition}

\proof
Let $f:G\to \varprojlim G/G_E$ send $g$ to $([g]_E)$. Note $f$ is a homomorphism. Suppose $f(g)=f(h)$ for some $g,h\in G$. Then $gh^{-1}\in G/G_E$ for each basic $E$. Now for each entourage $D$ of $X$ there is a basic entourage $E$ of $X$ so that the orbits of the action of $G_E$ on $X$ are $D$-bounded. Therefore for each $x\in X$, $(x,gh^{-1}x)\in D$ for each entourage $D$ of $X$. Since $X$ is Hausdorff $x=g^{-1}hx$ and since the action is faithful $g=h$. Therefore $f$ is injective and $G$ is isomorphic to a subgroup of $\varprojlim G/G_E$.

Now suppose the projection induced from the action of $G$ on $X$ has complete fibers. To see that $f$ is surjective and therefore $X$ is isomorphic to $\varprojlim G/G_E$, suppose $([g_E]_E)\in \varprojlim G/G_E$. Let $x\in X$ and set $A_E=\{g_Fx:F\subset E\}$. Note that each $A_E$ is contained in the fiber of $[x]\in X/G$. To see the $\{A_E\}$ is Cauchy, given a basic entourage $E$ of $X$ choose a basic entourage $F$ so that the orbits of the action of $G_F$ on $X$ are $E$-bounded. Suppose $g_Hx,g_Kx\in A_F$. Now $g_K^{-1}g_H\in G_F$ so $(g_K^{-1}g_Hx,x)\in E$. Since $E$ is invariant $(g_Hx,g_Kx)\in E$. Therefore $\{A_E\}$ is Cauchy so there is a limit point $gx$. Now, given an entourage $E$ of $X$ there is an entourage $F$ of $X$ so that $A_F\subset B(gx,E)$. Then $A_{F\cap E}\subset B(gx,E)$ so $(g_{F\cap E}x,gx)\in E$. Since $E$ is invariant $(g^{-1}g_{F\cap E}x,x)\in E$ so $g^{-1}g_{F\cap E}\in G_E$. Then $[g_{F\cap E}]_E=[g]_E$. But $[g_{F\cap E}]_E=[g_E]_E$ so $f(g)=([g_E]_E)$.
\endproof

\begin{proposition} \label{X=limX/G_E}
Suppose a group $G$ acts uniformly equicontinuously on a uniform space $X$. Suppose $X$ is Hausdorff and the action has small scale bounded orbits. Consider the inverse system $\{X/G_E,\phi_{FE}\}$ given by the basis as above. Then $X$ embeds in $\varprojlim X/G_E$. If the projection $X\to X/G$ induced from the action of $G$ on $X$ has complete fibers then $X$ is uniformly equivalent to $\varprojlim X/G_E$.
\end{proposition}

\proof
Let $f:X\to \varprojlim X/G_E$ send $x$ to $([x]_E)$. First, to see that $f$ is injective, suppose $f(x)=f(y)$ for some $x,y\in X$. Now for each entourage $D$ of $X$ there is a basic entourage $E$ so that the orbits of the action of $G_E$ on $X$ are $D$ bounded. Since $f(x)=f(y)$, there is a $g_E\in G_E$ with $x=g_Ey$. Then $(x,y)=(g_Ey,y)\in D$. Since $(x,y)\in D$ for each entourage $D$ of $X$ and $X$ is Hausdorff, $x=y$. To see that $f$ is a uniform embedding, first note it is uniformly continuous since the projections $X\to X/G_E$ are uniformly continuous. Suppose $D$ is an entourage of $X$. Choose $E$ so that $E^2\subset D$ and choose a basic entourage $F$ so that the orbits of the action of $G_F$ on $X$ are $E$-bounded. Suppose $x,y\in X$ with $(f(x),f(y))\in \pi_F^{-1}(p_F(E))$ where $\pi_F$ is the projection from $\varprojlim X/G_E$ to $X/G_F$ and $p_F$ is the projection associated with the action of $G_F$ on $X$. Now $(x,g_Fy)\in E$ for some $g_F\in G_F$. Also $(g_Fy,y)\in E$ so $(x,y)\in E^2\subset D$ and $(f(x),f(y))\in f(D)$.

Finally, to see that if $f$ is surjective and therefore a uniform equivalence if the projection induced from the action of $G$ on $X$ has complete fibers, suppose $([x_E]_E)\in \varprojlim X/G_E$. For each basic entourage $E$ of $X$ let $A_E$ be the orbit of $x_E$ under the action of $G_E$. Now each $x_E$ gets sent to the same equivalence class in $X/G$, say $[x]$, and each $A_E$ is contained in the fiber $p^{-1}([x])$ where $p$ is the projection associated with the action of $G$ on $X$. Note $\{A_E\}$ is a Cauchy filter base since given an entourage $D$ of $X$ there is a basic entourage $E$ of $X$ so that the orbits of the action of $G_E$ on $X$ are $D$-bounded. Therefore there is a limit $y$ of $\{A_E\}$. Let us see that $f(y)=([x_E]_E)$. Given $E$, there is an $A_F\subset B(x,E)$. Then $A_{F\cap E}\subset B(y,E)$. Now $y=g x_{F\cap E}$ for some $g\in G$. Then $(x_{F\cap E},g x_{F\cap E})\in E$ so $g\in G_E$. Therefore $[y]_E=[x_{F\cap E}]_E$. But $[x_E]_E=[x_{F\cap E}]$.
\endproof

\begin{theorem}
Suppose a group $G$ acts faithfully and uniformly equicontinuously on a metrizable uniform space $X$. Then the following are equivalent.
\begin{itemize}
\item[1.] The action has small scale bounded orbits and the projection $p:X\to X/G$ has complete fibers.
\item[2.] There is a Mittag-Leffler inverse sequence of groups $\{G_i,\psi_i\}$ with trivial $\ilim^1$, an inverse sequence of uniform spaces $\{X_i, \phi_i\}$, and compatible actions of $G_i$ on $X_i$ that are uniformly properly discontinuous and neutral with $G$ isomorphic to $\varprojlim G_i$, $X$ uniformly equivalent to $\varprojlim X_i$, and the action of $G$ on $X$ equivalent to the induced action of $\varprojlim G_i$ on $\varprojlim X_i$.
\end{itemize}
\end{theorem}

\proof
1. $\implies$ 2. Since $X$ is metrizable it is Hausdorff. Then by \ref{BasisOfPropDiscActions} and the paragraph succeeding it there are inverse systems $\{G/G_E, \psi_{FE}\}$ and $\{X/G_E,\phi_{FE}\}$ and compatible actions of $G/G_E$ on $X/G_E$ that are uniformly properly discontinuous. Since $X$ is metrizable it has a countable basis of entourages so these systems can be realized as sequences. Notice the actions are neutral since the action of $G$ on $X$ is uniformly equicontinuous. By \ref{G=limG/G_E} $G$ is isomorphic to $\varprojlim G_\alpha$. By \ref{X=limX/G_E} $X$ is uniformly equivalent to $\varprojlim G/G_E$. The inverse sequence $\{G/G_E, \psi_{FE}\}$ is Mittag-Leffler and has trivial $\ilim^1$ since in fact each $\psi_{FE}$ is surjective. Indeed, given $F\subset E$ and $[g]_E\in G/G_E$, $[g]_F\in G/G_F$ gets mapped to $[g]_E$. Notice the action of $G$ on $X$ is equivalent to the induced action of $\varprojlim G/G_E$ on $\varprojlim X/G_E$ since given $x\in X$ and $g\in G$, $([gx]_E)=([g]_E)([x]_E)$.

2. $\implies$ 1. By \ref{InvLimOfSmScBddOrbits} the action of $G$ on $X$ has small scale bounded orbits. According to the proof of \ref{InverseLimitOfRegUCMs} $p$ has complete fibers.
\endproof

\end{document}